\pdfoutput=1
\documentclass[a4paper]{article}

\usepackage{lipsum}
\usepackage{amsfonts}
\usepackage{graphicx}
\usepackage{epstopdf}
\usepackage{algorithmic}
\usepackage{amsmath}
\usepackage{amsopn}
\usepackage[]{cite} 
\usepackage{graphicx}
\usepackage{subcaption}
\usepackage{tikz}
\usepackage{url} 
\usepackage[hidelinks]{hyperref}
\usepackage{cleveref}
\ifpdf
  \DeclareGraphicsExtensions{.eps,.pdf,.png,.jpg}
\else
  \DeclareGraphicsExtensions{.eps}
\fi

\newcommand{\TheTitle}{Matrix-Free Evaluation Strategies for Continuous and Discontinuous Galerkin Discretizations on Unstructured Tetrahedral Grids}

\title{{\TheTitle}\thanks{The research presented in this paper was partly funded by
the German Ministry of Research, Technology and Space through project “PDExa: Optimized
software methods for solving partial differential equations on exascale supercomputers”, grant agreement no. 16ME0637K, and the European Union – NextGenerationEU, and by BREATHE, an ERC-2020-ADG Project, Grant Agreement ID 101021526.}}
\author{Dominik Still\thanks{Institute for Computational Mechanics,
TUM School of Engineering and Design, Technical University of Munich
  (\texttt{dominik.still@tum.de}, \texttt{wolfgang.a.wall@tum.de}).}
\and Niklas Fehn\thanks{High-Performance Scientific Computing
University of Augsburg - Faculty of Mathematics, Natural Sciences, and Materials Engineering 
  (\texttt{niklas.fehn@uni-a.de}).}
\and Wolfgang A. Wall\footnotemark[2]
\and Martin Kronbichler\thanks{
Ruhr University Bochum - Faculty of Mathematics 
  (\texttt{martin.kronbichler@rub.de}).}}

\newcommand{\software}{\texttt}

\newtheorem{remark}{Remark}[section]

\begin{document}

\maketitle

\begin{abstract}
This study presents novel strategies for improving the node-level performance of matrix-free evaluation of continuous and discontinuous Galerkin spatial discretizations on unstructured tetrahedral grids. In our approach the underlying integrals of a generic finite-element operator are computed cell-by-cell through numerical quadrature using tabulated dense local matrices of shape functions, achieving high throughput for low to moderate-order polynomial degrees. By employing dense matrix-matrix products instead of matrix-vector products for the cell-wise interpolation, the method reaches over $60\%$ of peak performance. The optimization strategies exploit explicit data parallelism to enhance computational efficiency, complemented by a hierarchical mesh reordering algorithm that improves data locality. The matrix-free implementation achieves up to a $6\times$~speedup compared to a global sparse matrix-based approach at a polynomial degree of three.

The effectiveness of the method is demonstrated through numerical experiments on the Poisson and Navier--Stokes equations. The Poisson operator is preconditioned by a hybrid multigrid scheme that combines auxiliary continuous finite-element spaces, polynomial and geometric coarsening where possible while employing algebraic multigrid on the coarse mesh. Within the preconditioner, the implementation transitions between the matrix-free and matrix-based strategies for optimal efficiency. Finally, we analyze the strong scaling behavior of the Poisson and Helmholtz operators, demonstrating the method's potential to solve large real-world problems.
\end{abstract}


\section{Introduction}
\label{sec:intro}

\begin{sloppypar}
For many applications of the finite element method, the solution of linear systems of equations is the most time-consuming algorithmic component. Large-scale problems are usually addressed by iterative solvers with suitable preconditioners for optimal computational complexity, where matrix-vector products are of crucial importance. In order to optimize the performance of these computational kernels, matrix-free operator evaluation has been established as an attractive approach for both continuous (CG) \cite{Kronbichler2012, Böhm2024, Moxey2020} and discontinuous Galerkin (DG) \cite{Kronbichler2019, Bergbauer2024} discretizations. By evaluating the local integrals of a finite element weak form on a global vector through a loop over cells, it is not necessary to construct a global sparse matrix. However, the optimal evaluation strategy for the local contributions across different polynomial degrees and element types is not clear a priori~\cite{Cantwell2010,Cantwell2011}, where the present work aims to provide additional insights.

To accelerate the evaluation of the local integrals for quadrature formulas with tensor-product structure, sum-factorization techniques are a common approach~\cite{Karniadakis1999, Kronbichler2012, Kronbichler2019}, providing a computational complexity of $\mathcal{O}(p^{d+1})$ per cell for polynomial degree $p$ in $d$ dimensions, improving over the complexity with dense interpolation matrices of $\mathcal{O}(p^{2d})$. For hypercube elements, applying sum-factorization is straightforward. However, hexahedral mesh generation on arbitrary domains remains an open question \cite{Pietroni2022}, while robust algorithms exist to generate tetrahedral meshes \cite{Shewchuk1998, Schneider2022}.

Moxey et al.~\cite{Moxey2020} present state-of-the-art implementations of sum-factorization algorithms on simplex elements for low to high polynomial degrees. Yet, early studies suggest that local dense interpolation matrices are more efficient for low to moderate orders~\cite{Cantwell2010,Cantwell2011}. This raises the question of the optimal evaluation strategy for different polynomial degrees on tetrahedral grids.
\end{sloppypar}
For low order polynomials on block structured tetrahedral grids, information embedded in the grid can be leveraged to create efficient computational kernels for the matrix-vector products in a matrix-free setting with few local stencils, avoiding the storage of all matrix entries~\cite{Kohl2023,Böhm2024}. A code generation framework utilizing several optimizations is developed in these studies, including inter-element vectorization, under-integration, moving of loop invariants, and specific loop patterns. However, these methods rely on exploiting the block structure of the grid and come with challenges for the mesh generation process, akin to hexahedral meshing.

Fully unstructured grids offer greater flexibility for complex geometries, but partially different optimization strategies have to be employed. For instance~\software{PyFR}~\cite{Witherden2013} expresses element-local computations as matrix-matrix products using pre-tabulated basis functions, leading to dense matrix-matrix multiplications for flux-reconstruction and solution updates. The operations are done concurrently for all points within an element, with multiple elements batched together for efficiency. The methods are also suitable for GPUs~\cite{Witherden2013,Akkurt2021}, which offer broad support for parallel computations~\cite{libceed}, including applications to DG methods~\cite{Kloeckner2009}. These systems provide data level parallelism through the Single Instruction, Multiple Data (SIMD) paradigm. On CPUs, various vectorization strategies are possible \cite{Sun2020}. Inter-element vectorization~\cite{Kronbichler2012} is applied in this work, as in \cite{Moxey2020} for simplex elements.

Optimizations can likewise be applied to the construction of element local matrices. For affine mappings, the geometry information can be pulled out of the weak form integral, producing a tensor that contains the contributions evaluated on the reference element only \cite{Kirby2005}. By reformulating the integral into a series of tensor contractions between the tensor and the cell specific geometry, the arithmetic complexity of the cell assembly and, in turn, the local operator evaluation, can be reduced. Additional gains are possible by exploiting the sparsity structure of the tensor. However, the decomposition is nontrivial to generalize to higher-order mappings and it is also not clear whether possible arithmetic gains extend to arbitrary variational forms.

This work contributes by analyzing the matrix-free finite element operator evaluation based on quadrature with dense interpolation matrices for low to moderate order polynomials and comparing its performance to that of a global sparse matrix approach.
By identifying generic building blocks~\cite{Kronbichler2012,Knepley2013,libceed}, the operator is flexible on the type of equation: Local interpolation matrices are populated with values and derivatives of basis functions at quadrature points, where the weak form is implemented.
Our evaluation approach focuses on node-level performance optimizations for continuous and discontinuous Galerkin methods on unstructured tetrahedral grids. We achieve performance gains by leveraging explicit data parallelism with SIMD techniques. Matrix-matrix products replace element-wise matrix-vector products, aligning the approach with the compute to data transfer ratio of modern CPU systems. Further, we introduce a hierarchical mesh reordering, improving the data locality on unstructured grids produced by typical mesh generators. The implementation is based on \software{deal.II}~\cite{dealii}, which has been shown to scale well on massively parallel, distributed-memory compute systems~\cite{Bangerth2012}. The matrix-free methodology is designed to integrate efficiently with hybrid multigrid (MG) preconditioners, using matrix-free operator evaluations and level transfers~\cite{Fehn2019,Munch2022}, enabling efficient solvers for Poisson and Helmholtz-type problems, which are fundamental in the context of the Navier--Stokes equations.

This paper is organized as follows: \Cref{sec:matrix_free} details the optimization strategies for the matrix-free operator, while \Cref{sec:numerical_experiments} shows its properties in comparison to global sparse matrices. \Cref{sec:application} presents experimental results of the hybrid-MG-preconditioned Poisson operator and the Navier--Stokes equations. \Cref{sec:conclusions} concludes the study with a summary of key findings.

\section{Matrix-Free Operator Evaluation} \label{sec:matrix_free}
In this section, matrix-free strategies for continuous and discontinuous Galerkin discretizations of the Poisson operator are presented. Based on memory access and arithmetic requirements, a series of optimizations are performed.

\subsection{Poisson Equation Formulations}
Consider the Poisson equation on an arbitrary domain~$\Omega \subset \mathbb{R}^d$
\begin{equation}\label{eq:poisson}
  - \Delta u = f \text{ in } \Omega, \text{ with } u = u_D \text{ on } \Gamma = \partial \Omega \; .
\end{equation}
The weak form seeks the field $u$ that fulfills, for all test functions~$v$, the relation
\begin{equation}\label{eq:weak_form_CG}
  a_{\mathrm{CG}}(v,u) = {(\nabla u , \nabla v)}_{\Omega} = b(v) .
\end{equation}
On a discretized domain~$\Omega_h$, we can solve the problem on either the continuous Galerkin finite element subspace of~$H^1(\Omega)$ or with the discontinuous Galerkin method on a finite-dimensional subspace of~$L^2(\Omega)$.

For discontinuous Galerkin methods, the continuity of the elements is enforced by the flux over the element faces. Considering the symmetric interior penalty Galerkin~(SIPG) method~\cite{Hesthaven2007}, the bilinear form reads
\begin{equation} \label{eq:weak_form_DG}
    a_{\mathrm{DG}}(u,v) = {(\nabla u , \nabla v)}_{\Omega_h} - ( \{\!\{ \partial_n v \}\!\}, [\![u]\!])_{ {\Gamma}_{h}^{\text{int}}} - (  [\![v]\!], \{\!\{\partial_n u\}\!\} - \tau [\![u]\!] )_{ {\Gamma}_{h}^{\text{int}}} \; ,
\end{equation}
where~$[\![\cdot]\!] = (\cdot)^- - (\cdot)^+$ is the jump operator and~$\{\!\{\cdot\}\!\} = ((\cdot)^- + (\cdot)^+)/2$ the average operator at interior faces in~$\Gamma_{h}^{\text{int}}$ and $\tau$ is the interior penalty parameter. At the boundary, suitable modifications are made to $a_{\mathrm{DG}}$ and $b(v)$ \cite{Hesthaven2007}.

\subsection{Matrix-Free Operator}

Inserting ansatz and test functions and integrating leads to a linear system of equations of the form~$\textbf{A} \textbf{x} = \textbf{b}$.
To solve this system, iterative solvers are typically used, which apply the action of the system matrix~$\textbf{A}$ on some vector~$\textbf{u}$ as
\begin{equation}
    \textbf{y} = \textbf{A}\textbf{u} \; .
\end{equation}

Instead of assembling the matrix, we use a matrix-free approach~\cite{Kronbichler2012}, denoted here as the evaluation of an abstract operation~$\textbf{A}(\textbf{u})$. The implementation of this operator traverses all cells and adds the local contribution of a cell into a global result vector.
In the matrix-free setting and considering a continuous Galerkin discretization, the operator evaluation may be written as
\begin{equation} \label{eq:matrix_free_cg}
    \textbf{A}(\textbf{u}) = \sum_{e = 1}^{N_{el}} \textbf{G}^{\text{T}}_e \textbf{E}_e^{\text{T}} \textbf{D}_e \textbf{E}_e \textbf{G}_e \textbf{u} \; .
\end{equation}

\begin{figure}
    \centering
    \includegraphics[width=\textwidth]{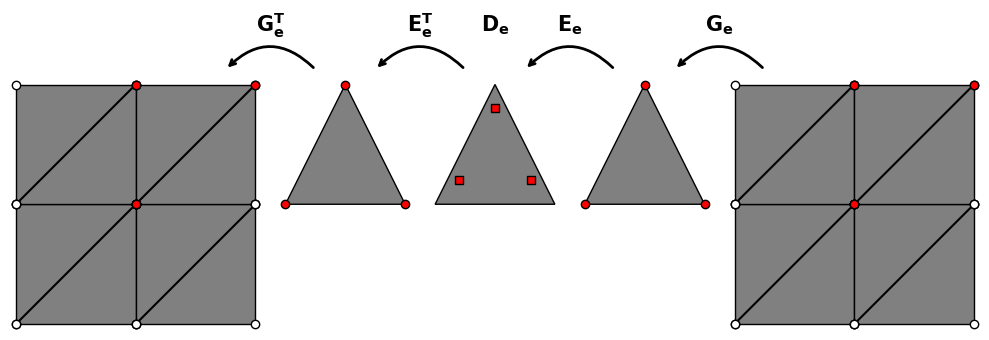}
    \caption{Illustration of main steps of matrix-free operator evaluation on a 2D simplicial mesh~(adopted from~\cite{Fehn2019}). First, the operator $\textbf{G}_e$ selects the relevant DoF values. Then, the matrix $\textbf{E}_e$ evaluates the gradients at the quadrature points. The operator $\textbf{D}_e$ first transforms the gradient with the inverse of the Jacobian from the reference cell to the real cell, then applies the quadrature weights and the determinant of the Jacobian and then transforms the gradient back to the reference cell. The matrix $\textbf{E}_e^T$ interpolates back to the DoF values and the operator $\textbf{G}^{\text{T}}_e$ scatters the DoF values back into the global vector.}
    \label{fig:matrix_free_operator}
\end{figure}

Figure \ref{fig:matrix_free_operator} illustrates this process. Starting from the right, the operator~$\textbf{G}_e$ selects the entries of the global vector~$\textbf{u}$ that are associated with the element $e$. These entries are referred to as the locally relevant degrees of freedom (DoF) values, denoted as $\mathbf{u}_e$.

The elementwise coefficients are then used to evaluate the gradients of the solution field $\nabla u_h(\textbf{x}) = \sum_j \nabla \phi_j(\textbf{x}) u_{e,j}$ at quadrature points $\hat{\textbf{x}}_i$, $i=1,\ldots, n_q$. Using the transformation to reference coordinates $\textbf{x}\mapsto \hat{\textbf{x}}$, the gradients are evaluated as $\nabla \phi_j(\hat{\textbf{x}}(\textbf{x})) = \textbf{J}_e(\hat{\textbf{x}})^{-\text T}\hat{\nabla} \hat{\phi}_j(\hat{\textbf{x}})$ with the Jacobian matrix of the map from reference to real coordinates denoted by $\textbf{J}$ and $\hat{\nabla} = \left(\hat{\partial}_1, \ldots, \hat{\partial}_d\right)$ being the gradient with partial derivatives with respect to the coordinates on the reference simplex. Hence, we represent an evaluation step by the matrix~$\textbf{E}_e$ as
\begin{equation} \label{eq:interpolation_matrix}
    (\textbf{E}_e)_{di+k,j} = \hat{\partial}_k \hat{\phi}_j(\textbf{$\hat{\textbf{x}}$}_i) \; \text{with} \; k=1,...,d \;.
\end{equation}
Here,~$\hat{\phi}_j$ is the $j^{\text{th}}$ shape function, with $j = 1,..., n_{\text{DoF}}$ for all  $n_{\text{DoF}}$ shape functions (or DoFs) on the cell. 
By definition, the matrix~$\textbf{E}_e$ is the same on every cell as its action computes the derivatives on the unit cell for all quadrature points from all local DoF values, seen by its size of $n_{\text{q}} \cdot d \times n_{\text{DoF}}$. Note that the multiplication with the Jacobian $\mathbf{J}$ is shifted to the next step.

The operator~$\textbf{D}_e$ acts at each quadrature point by first applying the inverse transpose Jacobian $\textbf{J}_e(\hat{\textbf{x}})^{-\text T}$ to obtain the gradient in real space, then multiplying by the quadrature weight and the determinant of the Jacobian, and finally applying the inverse Jacobian to transform the gradient back.

The transpose~$\textbf{E}_e^{\text{T}}$ operator represents multiplication by derivatives of test functions and summation over all quadrature points. Finally, the operator~$\textbf{G}^{\text{T}}_e$ scatters the local DoF values back into the global DoF vector. 

For discontinuous spaces of shape functions, an additional loop over all faces is required, in accordance to~\cite{Kronbichler2019} it is written as a separate loop,
\begin{equation} \label{eq:matrix_free_dg}
    \textbf{A}(\textbf{u}) = \left(\sum_{e = 1}^{N_{\mathrm{el}}} \textbf{G}^{\text{T}}_e \textbf{E}_e^{\text{T}} \textbf{D}_e \textbf{E}_e \textbf{G}_e+ \sum_{f = 1}^{N_{\mathrm{faces}}} \textbf{G}^{\text{T}}_f \textbf{E}_f^{\text{T}} \textbf{D}_f \textbf{E}_f \textbf{G}_f \right) \textbf{u}.
\end{equation}
In our code, the cell and face loops are interleaved to exploit data locality in the access to the vectors $\textbf{u}$ and $\textbf{y}$, i.e.~to reduce the amount of data transferred from main memory.

\subsection{Memory Transfer and Arithmetic Operations}
We aim to optimize the $\text{throughput} = \frac{\text{number of unknowns [$n_{\mathrm{DoF}}$]}}{\text{computational time [$\mathrm{s}$]}}$.
Considering the roofline model~\cite{Williams2009} as the performance-relevant model, we specify in this section the data transfer from main memory and the number of arithmetic operations. These estimates assume that all data associated with the global vectors and auxiliary structures, such as global indices, must be loaded from main memory, while optimal cache sizes and policies prevent further memory accesses for intermediate evaluation data and for repeatedly accessed entries of the vectors $\textbf{u}$ and $\textbf{y}$. Furthermore, boundary effects are neglected. We show results in three dimensions.

\subsubsection*{Continuous Galerkin} For the continuous discretization, we first establish a relation between the number of DoFs per cell and the global number of DoFs. For a~$d$-dimensional simplex element, the number of DoFs is $n_{\text{DoFs/cell}} = \frac{(p + d)!}{d!p!}$. We define the number of unique DoFs per cell as the ratio of the total number of DoFs to the total number of cells, which is the average number of distinct DoFs of each mesh element. It depends on the connectivity of the grid. Subdividing a unit cube into~5 tetrahedra~\cite{Fominykh2015,Benzley1995}, while accounting for the contributions of shared DoFs at the vertices, edges, faces, and the DoFs within the tetrahedron, leads to the number of unique DoFs per cell~$n_{\text{unique DoFs/cell}} = 0.2, ~1.4,~4.6$ for polynomial degrees of~$p=1,~2,~3$, respectively.
Bergbauer et al.~\cite{Bergbauer2024} estimate the memory transfer and arithmetic complexity for a single matrix-vector product. We extend \cite[Equations (5.6) and (5.7)]{Bergbauer2024} to tetrahedral elements with the aforementioned quantities, leading to a memory transfer $d_{\text{cell}}^{\text{CG}}$ byte per cell in double precision as
\begin{align} \label{eq:mem_trans_cg}
 \hspace{0,5cm} d_{\text{cell}}^{\text{CG}} =  \underbrace{\frac{14}{n_{\text{lanes}}}}_{\mathbf{D}_e \; \text{(geometry)}} + \underbrace{ 4 \cdot n_{\text{DoFs/cell}}}_{\mathbf{G}_e\; \text{(indices)}} +  \underbrace{3 \cdot 8 \cdot n_{\text{unique DoFs/cell}}}_{\textbf{u}, \textbf{y} \; \text{(read/write DoFs)}} \; .
\end{align}

The 14 bytes consist of indices and pointers to geometry information ($\textbf{D}_e$). Like in~\cite{Bergbauer2024}, we assume a SIMD evaluation strategy working on multiple cells in different lanes within the same instruction. The memory cost is amortized over the number of SIMD lanes~$n_{\text{lanes}}$ by setting pointers and indices concurrently for multiple cells. For each DoF on a cell, the global index has to be read ($\textbf{G}_e$, 4 bytes). A total of $3 \times 8$ bytes are from reading the DoF value once ($\mathbf{u}$) and reading (for ownership) and writing to the global vector ($\mathbf{y}$) in double precision (8 bytes each). If two cells share a DoF, the entries are assumed to be read only once, so each cell needs to read only its distinct number of DoFs $n_{\text{unique DoFs/cell}}$.

This model assumes that the interpolation matrix $\textbf{E}_e$ is fetched from cache. The Jacobian matrices $\mathbf{J}_e$ are precomputed, with identical matrices being identified so that only distinct ones need to be stored, thereby compressing the geometry data and optimizing memory usage \cite{Kronbichler2012}. In the case of affine mappings, a single Jacobian per cell suffices and is stored accordingly. As the Jacobians can potentially be reused multiple times, they are assumed to be available from cache.

Equation~\eqref{eq:matrix_free_cg} is used to determine the arithmetic costs. For each cell, the product of the matrix~$\textbf{E}_e$ and its transpose with the vector~$\textbf{u}_e$, along with the operations on the quadrature points, requires $w_{\text{cell}}^{\text{CG}}$ arithmetic operations,
\begin{align} \label{eq:arithm_cg}
  \hspace{0,5cm} w_{\text{cell}}^{\text{CG}}  =  \underbrace{2}_{\textbf{E}_e/\textbf{E}_e^{\text{T}}} \cdot \; \underbrace{3 \; n_{q} \;n_{\text{DoFs/cell}} \cdot 2}_{\textbf{E}_e \textbf{u}_e}  + \underbrace{(15 \cdot 2 + 3) n_{q}}_{\textbf{D}_e \; \text{(Jacobian matrix and JxW)}} + \underbrace{n_{\text{DoFs/cell}}}_{\textbf{y}}  \; .
\end{align}
The leading factor of $2$ is from applying the interpolation matrix $\mathbf{E}_e$ and its transpose $\mathbf{E}_e^\text{T}$. The matrix-vector product $\mathbf{E}_e \mathbf{u}_e$ has the cost $3 n_q$ (rows in the matrix) multiplied by $n_{\text{DoFs/cell}}$ (length of the vector), the second factor of $2$ is from the multiply-add operations. For every quadrature point, the $3 \times 3$ Jacobian and its transpose are applied to the gradient vector (15 operations for each matrix-vector product) and the vector is scaled by JxW, a precomputed factor storing the product of quadrature weight and determinant of the Jacobian (operation $\mathbf{D}_e$). The number of quadrature points $n_{q}$ is given by the quadrature rule, adopted from \cite{Witherden2014} and \cite{Xiao2010}, chosen to exactly integrate polynomials of degree $p+1$ for affine mappings. The last term stems from summation into the global vector $\textbf{y}$.

\subsubsection*{Discontinuous Galerkin} The memory transfer per cell is given as
\begin{equation} \label{eq:mem_trans_dg}
\hspace{0,5cm} d_{\text{cell}}^{\text{DG}} =  \underbrace{\frac{14+4}{n_{\text{lanes}}}}_{\mathbf{D}_e, \; \mathbf{G}_e} + \underbrace{4}_{n_{\text{faces}}} \underbrace{\left(4+ \frac{18}{n_{\text{lanes}}}\right)}_{\mathbf{G}_f, \mathbf{D}_f \; \text{(face geometry)}}  \; 
+ \underbrace{3 \cdot 8 \cdot n_{\text{DoFs/cell}}}_{\textbf{u}, \textbf{y}\; \text{(read/write DoFs)}} \; .
\end{equation}
For the cell geometry, the same data as in the continuous case is accessed ($\mathbf{D}_e$). Since DoFs are uniquely associated with individual cells, one DoF index per cell is loaded ($\mathbf{G}_e$), while the remaining DoF indices can be reconstructed from the first DoF index~\cite{Kronbichler2019}. For each of the 4 $n_{\text{faces}}$ of a tetrahedron, 4 bytes of indices are loaded and analog to the cell term, 18 bytes of face geometry information ($\mathbf{D}_f$) are needed. Data access to the global vectors is identical to the continuous case.

Regarding the number of arithmetic operations, additional face terms have to be evaluated for the SIPG method. Assuming an affine mapping, we obtain the number of additional operations as
\begin{equation} \label{eq:arithm_dg}
    \begin{split}
 \hspace{0,5cm} w_{\text{cell}}^{\text{DG}}  =  w_{\text{cell}}^{\text{CG}} + 2 \cdot \underbrace{2}_{\textbf{E}_f/\textbf{E}_f^{\text{T}}} \cdot \; \underbrace{(1 + 3) \; n_{q_f}\; n_{\text{DoFs/cell}} \cdot 2}_{\textbf{E}_f \textbf{u}_f}  + \hspace{-0.5cm}\underbrace{43 n_{q_f}}_{\text{jump/average operations}} .
    \end{split}
\end{equation}
On the faces, two matrix-vector products have to be evaluated, one to interpolate the values and one to evaluate the gradients at the quadrature points, resulting in two matrices of size $n_{q_f} \times n_{\text{DoFs/cell}}$ and $3 n_{q_f} \times n_{\text{DoFs/cell}}$. Each face is shared by two neighboring cells, leading to a factor of two. Computing the jump and the averaged normal derivative takes 43 operations on each face quadrature point $n_{q_f}$.

\begin{figure}
    \centering
   \begin{subfigure}{0.5\textwidth}
  \centering
  \includegraphics[width=0.99\linewidth]{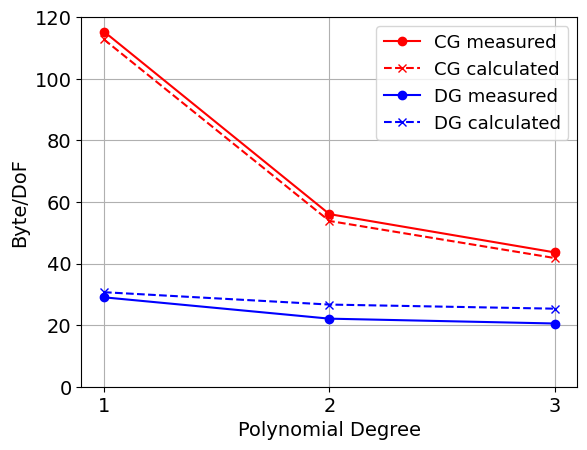}
  \caption{Memory transfer}
  \label{fig:mem}
\end{subfigure}%
\begin{subfigure}{.5\textwidth}
  \centering
  \includegraphics[width=0.99\linewidth]{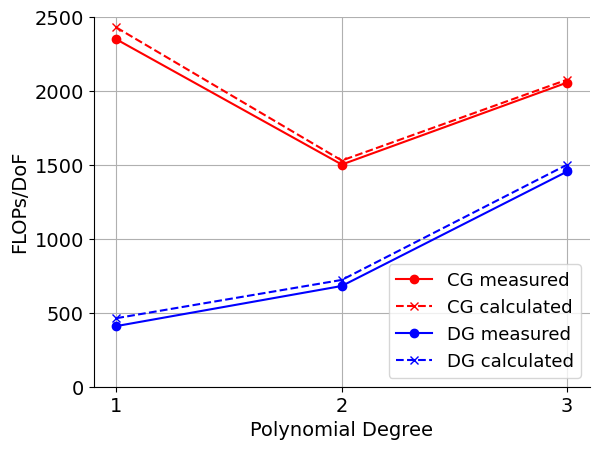}
  \caption{Number of arithmetic operations}
  \label{fig:ops}
\end{subfigure}
    \caption{Comparison of theoretical and measured performance metrics on a simplex grid in 3D.}
    \label{fig:arith_mem}
\end{figure}

The predictions are in good agreement with measurements from hardware performance counters evaluated through the tool~\software{LIKWID}~\cite{likwid}, as shown in Figure \ref{fig:arith_mem}. The measurements are done on a cube geometry discretized with 9.3 million cells. Note that the reduced number of arithmetic operations is due to the reduced work near the boundary, which is omitted in the prediction.

\begin{remark}
    Note that the arithmetic complexity of evaluating the local operator using tensor-based representations for affine geometries according to~\cite{Kirby2005} mentioned in the introduction would have a lower constant, $18 (n_{\text{DoFs/cell}}^2 + n_{\text{DoFs/cell}})$ compared to the cost in Equation~\eqref{eq:arithm_cg}. Further gains would be possible due to the matrix structure, see also the sub-structuring exploited by \cite{Uphoff2016}. However, since the generalization to curved geometries, variable coefficients, other variational forms or face integrals is not clear, we do not consider alternatives to the quadrature approach in this work.
\end{remark}

\subsection{Optimization of Arithmetic Work for Interpolation Operation} \label{sec:operator_optimizations}
\begin{sloppypar}
From the analysis of the arithmetic costs (Equations~\eqref{eq:arithm_cg} and~\eqref{eq:arithm_dg}), it can be deduced that the matrix-vector product~$\textbf{E}_e \textbf{u}_e$ dominates the arithmetic costs.
\end{sloppypar}
\subsubsection*{Reference Implementation} \label{sec:reference_implementation}
Looking at Equation~\eqref{eq:interpolation_matrix}, the shape function values and derivatives are evaluated at the quadrature points on the reference element. Since the interpolation operations represented by matrix $\textbf{E}_e$ are identical for all elements, see also~\cite{Kronbichler2012}, inter-element SIMD vectorization is employed to execute arithmetic instructions in Equations~\eqref{eq:matrix_free_cg} and~\eqref{eq:matrix_free_dg} for $n_{\text{lanes}}$ cells at once. Therefore, the cells (and faces if applicable) are grouped together in a setup phase in batches of size $n_{\text{lanes}}$.

\begin{figure}
    \centering
   \begin{subfigure}{0.33\textwidth}
  \centering
  \includegraphics[width=0.95\linewidth]{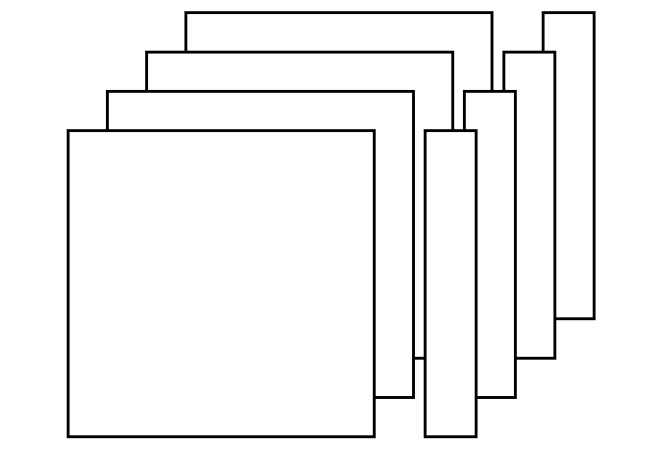}
  \caption{Reference}
  \label{fig:refernce_vectorisation}
\end{subfigure}%
\begin{subfigure}{.33\textwidth}
  \centering
  \includegraphics[width=0.95\linewidth]{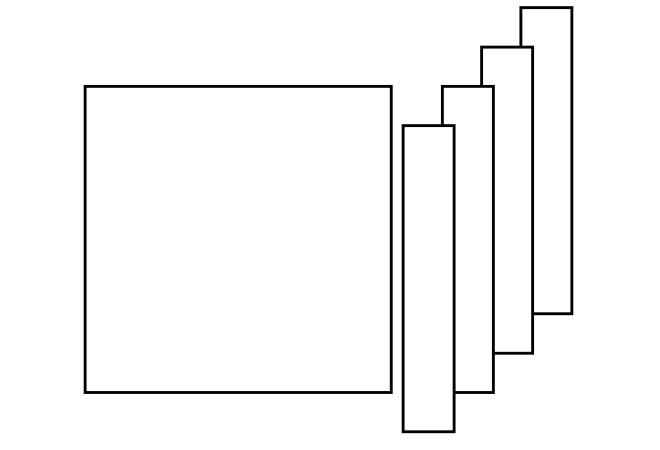}
  \caption{Data structures}
  \label{fig:data_vectorisation}
\end{subfigure}
\begin{subfigure}{.33\textwidth}
  \centering
  \includegraphics[width=0.95\linewidth]{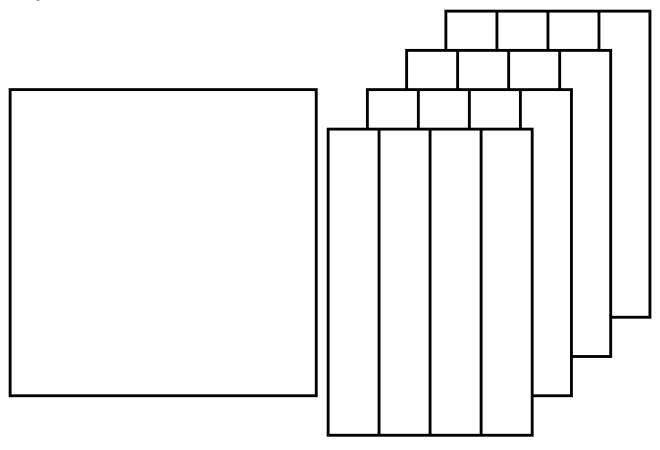}
  \caption{MM-Mult}
  \label{fig:mm_vectorisation}
\end{subfigure}
    \caption{In the reference state, $n_{\text{lanes}}$ copies (here shown with 4 SIMD lanes) of the interpolation matrix and $n_{\text{lanes}}$ DoF vectors are loaded (a), for broadcasting this changes to one copy of the interpolation matrix, while the number of vectors stays the same (b). In the last case, 4 DoF vectors are batched into a matrix, such one interpolation matrix and $4 \times n_{\text{lanes}}$ DoF vectors are loaded, increasing the reuse of matrix entries (c).}
    \label{fig:vectorisation}
\end{figure}

For each SIMD lane, the interpolation matrix and the DoF vector is loaded, as illustrated in Figure \ref{fig:refernce_vectorisation}. Notably, for tensor-product elements, loading a copy of the interpolation matrix for each SIMD lane was deemed efficient in \cite{Kronbichler2012, Kronbichler2019} to avoid the need of broadcast (scalar to SIMD vector) instructions, as only the one-dimensional interpolation matrix is used in the sum-factorization algorithm.

\begin{figure}
    \centering
   \begin{subfigure}{0.5\textwidth}
  \centering
  \includegraphics[width=0.99\linewidth]{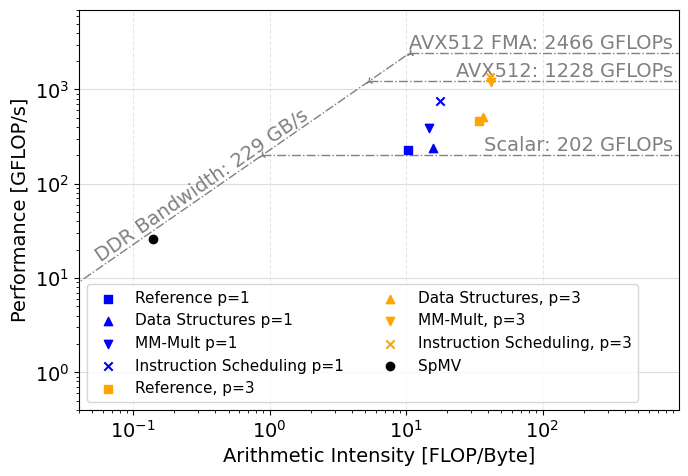}
  \caption{Roofline model}
  \label{fig:roofline_cg}
\end{subfigure}%
\begin{subfigure}{.5\textwidth}
  \centering
  \includegraphics[width=0.99\linewidth]{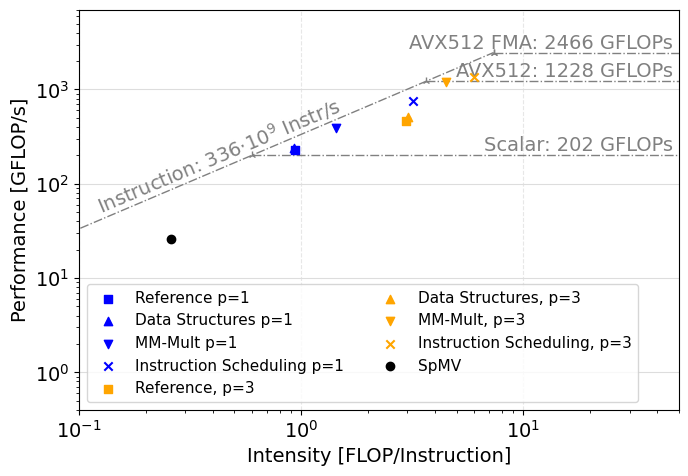}
  \caption{Instruction-roofline model}
  \label{fig:roofline_model_instructions}
\end{subfigure}
    \caption{Roofline model and instruction-roofline model of the continuous matrix-free operator at~$p=1$ and~$p=3$ compared to the sparse matrix-vector product (SpMV).}
    \label{fig:roofline}
\end{figure}

\begin{figure}
    \centering
   \begin{subfigure}{0.95\textwidth}
  \centering
  \includegraphics[width=0.6\linewidth]{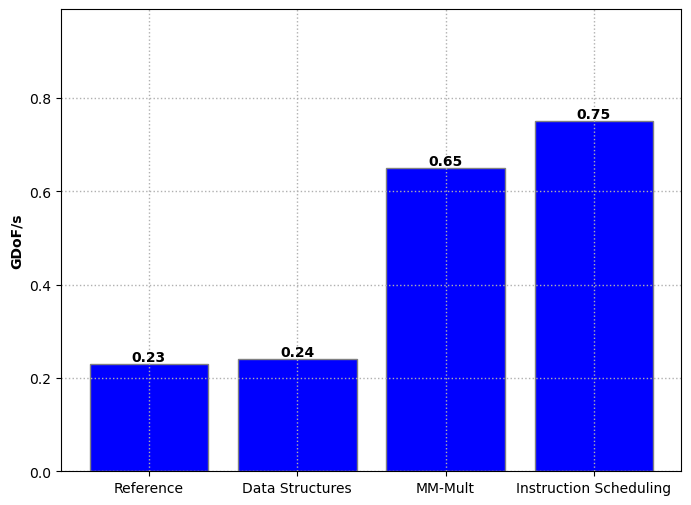}
  \label{fig:throughput_opt_CG}
\end{subfigure}
    \caption{Throughput under different optimizations at~$p=3$.}
    \label{fig:throughput_opt}
\end{figure}

\subsubsection*{Data Structures/Broadcasting} Given the size of the interpolation matrix, instead of loading ~$n_{\text{lanes}}$ copies of the matrix, a scalar storage is used and entries are broadcast to the full SIMD width of $n_{\text{lanes}}$ within the matrix-vector product with $\textbf{E}_e$ (see Figure \ref{fig:data_vectorisation}). This reduces the resident size of the interpolation matrix in the processor cache by a factor of~$n_{\text{lanes}}$ and possibly also cache bandwidth pressure. 
The SIMD vectorization is implemented using hardware intrinsics similar to \software{std::experimental::simd}~\cite{stdexperimentalsimd} and \software{XSIMD}~\cite{xsimd}. This enables compiler-supported vector arithmetic with near-optimal resource utilization while broadcasting. In the roofline model shown in Figure~\ref{fig:roofline_cg}, which exemplifies the continuous case for polynomial degrees $p=1$ and $p=3$, starting from the unoptimized operator (Reference), the arithmetic intensity increases (Data structures). Measurements are done on two Intel Xeon Gold~6230 sockets, see Table~\ref{tab:xeon_gold_platinum}, using \software{MPI} for distributed-memory computations. Each data point represents a complete matrix-vector product of the Poisson operator.

\begin{table}
\centering
\caption{Hardware properties of Intel Xeon Gold 6230 and Xeon Platinum 8360Y Processors}
\begin{tabular}{@{}lcc@{}}
\hline
\textbf{Property}               & \textbf{Xeon Gold 6230} & \textbf{Xeon Platinum 8360Y} \\
\hline
Microarchitecture               & Cascade Lake            & Ice Lake-SP                  \\
Base Clock Speed                & 2.10 GHz               & 2.40 GHz                     \\
Number of Cores (per socket)    & 2x20                   & 2x36                         \\
L1 Cache (per core)             & 32 kB                  & 48 kB                        \\
L2 Cache (per core)             & 1 MB                   & 1.25 MB                      \\
L3 Cache (shared)               & 24.75 MB               & 54 MB                        \\
Instruction Set                 & AVX-512                & AVX-512                      \\
\hline
\end{tabular}
\label{tab:xeon_gold_platinum}
\end{table}

\subsubsection*{Matrix-Matrix Product} The performance limit of the matrix-free algorithm working on a single element batch at a time arises not from loading the matrix entries from main memory but from streaming the data from the cache to the registers and arithmetic units. Hence, the optimization proposed in the previous subsection only marginally improves performance. To overcome this limitation, we switch from a matrix-vector product to a matrix-matrix product, which increases the reuse of matrix entries by register blocking. To convert the element-wise DoF vectors~$\textbf{u}_e$ into the matrix~$\textbf{U}_{e_1,...,e_n}$, two different strategies are explored in this work. In the first approach, the DoF vectors of four cells are arranged into a matrix, i.e.,\ we extend the cell batches. In Section~\ref{sec:reference_implementation} $n_{\text{lanes}}$ cells were batched to exploit the SIMD structure; now the cell batches are extended to hold $4 \cdot n_{\text{lanes}}$ cells, not only for the SIMD operations but also to be able to apply the matrix-matrix product, see Figure~\ref{fig:mm_vectorisation}. This means that in every SIMD lane, one matrix-matrix product is computed. The following holds in every SIMD lane for the matrix:
\begin{align}
    \textbf{U}_{e_1,...,e_4} = \left[\textbf{u}_{e_1}, \hdots, \textbf{u}_{e_4}\right],
\end{align}
where $e_1, \hdots, e_4$ are the indices of four cells.
The second approach addresses vector-valued problems. Instead of combining DoF vectors of different cells, the DoF vector of a single cell is reshaped into a matrix such that each column corresponds to one component of the vector-valued solution:   
\begin{equation}
    \textbf{U}_{1_e,2_e,...,d_e} = \left[{\textbf{u}_{{1}}}_{e},{\textbf{u}_{{2}}}_{e} ,\hdots, {\textbf{u}_{{d}}}_{e}\right]
\end{equation}
where~$\textbf{u}_{j_e}$ denotes the element-wise DoF vector associated with the~$j^{\mathrm{th}}$ component of the d-dimensional solution field. We refer to this as the components-batched strategy.

The size of both matrices is comparatively small, the interpolation operator~$\textbf{E}_e$ has dimension~$d \cdot n_q  \times n_\text{DoFs/cell}$, which, e.g., in 3D with polynomial degree~$p=3$, results in a matrix of size~$105 \times 20$. While dense matrix-matrix product algorithms like~\software{dgemm} reach~$85-90 \%$ of peak performance~\cite{Goto2008, Zee2015}, here custom compute kernels are used as the matrices fit into cache and we want to avoid function call and \software{dgemm} data rearrangement overheads. As the problem is reduced to small matrix-matrix products, optimization strategies similar to \software{LIBXSMM} \cite{Heinecke2016} are applied:
From plain C++ code, a compiler does not usually generate optimal code for the matrix-matrix product kernel. Thus, spatial blocking and ~$4 \times$ manual loop unrolling are employed, which enables the reuse of the entries of ~$\textbf{E}_e$. In the roofline model, the Poisson operator reaches approximately $50\%$ of the peak at $p=3$ (MM-Mult) for the cells-batched strategy, while the sparse matrix-vector (SpMV) product is memory-bound. 

Combined with the optimizations of the data structures, an increase in throughput can be seen in Figure~\ref{fig:throughput_opt}. The remaining gap to the arithmetic roofline is due to the other operators in Equation~\eqref{eq:matrix_free_cg}. Most prominently, the read/write phases involving the global vectors (operators $\textbf{G}_e$/$\textbf{G}_e^{\text{T}}$) contribute with instructions, but hardly any arithmetic work.

These phases can not effectively be overlapped with the matrix-matrix operations, resulting in a gap between the achieved and theoretical performance.
\begin{remark}
Strictly speaking, SIMD vectorization already transforms the operation into a matrix-matrix multiplication, as it aggregates $n_{\text{lanes}}$ DoF vectors into a matrix with $n_{\text{lanes}}$ columns. However, for the sake of clarity, we view it as performing a separate matrix-vector multiplication within each SIMD lane (see Figure~\ref{fig:data_vectorisation}).
\end{remark}

\subsubsection*{Instruction Scheduling} As noted in \cite{Moxey2020}, access to the global vector reduces the achievable performance. To mitigate this, we propose a more efficient approach compared to~\cite{Kronbichler2012} by eliminating the overhead associated with the implementation in the \software{deal.II} library~\cite{dealii}, which is designed to support multiple use cases. Specifically, we introduce a masked gathering function, which skips constrained DoFs, thereby reducing instruction overheads in the functions that selectively read and write the global DoF vector. Reducing the instruction overhead (Instruction Scheduling) pushes the operator at $p=3$ in the roofline model to around $60\%$ of the arithmetic peak. 
For linear shape functions, the operator is instruction-scheduling bound, as shown in Figure~\ref{fig:roofline_model_instructions}, where all data points of the linear case are near the instruction ceiling. In the instruction-roofline model, the number of total instructions per second is used instead of the memory bandwidth. This analysis assumes that the CPU can process four instructions per cycle, which is a heuristic number given the instruction mix in both phases of the algorithm and the support of a 4-wide decode of the CPU.

The optimizations mentioned above reduce the number of instructions needed for the computation. The matrix-matrix product reduces the number of multiplication kernel calls, and loop unrolling minimizes loop overhead within the kernel. Finally, instruction scheduling contributes to the improved performance at $p=1$, suggesting that our assumption of separate read/write phases and compute phases is correct. For cubic shape functions, the additional performance gains are explained by the fact that the compute phase was already compute-bound, while the read/write phase still could be improved, whereas in the linear case the instruction-mix is limiting, as the size of the interpolation matrix is not large enough to reach the computational limit. Note that the optimization of data structures/broadcasting does not change the number of instructions executed.

\begin{remark}
    We discussed the optimization of the matrix-vector product independent of the dimension of the problem. Looking at Equation~\eqref{eq:matrix_free_dg}, the same optimizations apply to the operations on the cell and the faces in the DG case.
\end{remark}

\subsection{Improving Cell Access by Hierarchical Grid Reordering} \label{sec:reordering}
Grid generators often produce meshes with irregular ordering, resulting in sub-optimal data locality due to scattered memory access. Hierarchically reordering the cells enhances data locality and improves cache utilization. This approach is analogous to matrix-based operators, where a sparsity pattern is created to efficiently store the matrix in terms of memory access patterns and data locality, as discussed in~\cite[Table 3]{Trotter2022}. 

In the matrix-free context, we introduce a hierarchical reordering of the mesh using a greedy algorithm designed to minimize bandwidth usage. Initially, neighboring cells are grouped together and the groups are ordered based on connectivity. Groups with higher connectivity are placed closer together in memory, while those with lower connectivity are positioned farther apart. This reordering process is applied recursively, ensuring that both local (within-group) and global (between-group) optimizations are achieved, mimicking the recursive structure of the Morton space filling curve~\cite{Morton1996}.

To illustrate this, consider an example where six cells are initially grouped into three groups. All three groups consist of two neighboring cells. The groups are then reordered by analyzing connectivity at the group level. If Group 1 has more connections to Group 3 than to Group 2, Group 3 is placed adjacent to Group 1, thus potentially improving efficiency.
The first loop optimizes the numbering of the cells within the group. By calling the function recursively, the ordering of the groups themselves is optimized, aiming to find an optimal ordering of groups to further improve efficiency.

\section{Throughput Studies for Operator Evaluation} \label{sec:numerical_experiments}
Having characterized the matrix-free continuous and discontinuous operators, we now analyze the performance of the operators in comparison to their matrix-based counterparts. Therefore, we consider the two strategies described above. Additionally, we analyze curvilinear elements as a borderline case with high data demand and we examine the case of a modified integration rule, where, in the matrix-free case, the number of computations can be reduced while maintaining the optimal convergence rate.

\subsection{Throughput of Cells-Batched and Components-Batched Strategies} \label{sec:thoughput}
\begin{figure}
    \centering
    \includegraphics[width=0.7\textwidth]{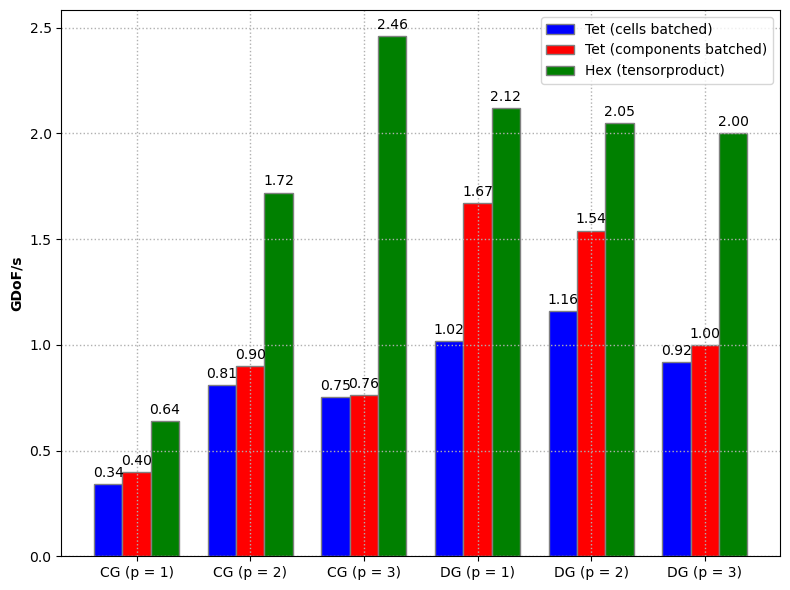}
    \caption{Throughput with different polynomial orders compared with hypercube elements.}
    \label{fig:throughput_gold}
\end{figure}

To evaluate the effect of different batching strategies for cells or vector components into a matrix, the throughput is analyzed under the assumption that the operator works in the saturated regime, i.e., most data has to be fetched from main memory. In the continuous case, the difference between the two strategies is minor, as seen in Figure~\ref{fig:throughput_gold}. However, the discontinuous case reveals distinct properties.

The cells-batched strategy requires reordering the faces on the fly and matching those with identical orientations, as they share the same interpolation matrices. This process is not always optimal, as finding four faces with the same orientation is not guaranteed, which can prevent the full utilization of optimized kernels. Grouping faces with the same orientation in a setup phase would be possible but at the price of data locality. The components-batched strategy avoids the issue, as the components on the face are assembled into the matrix, eliminating the need for face reordering. The components-batched strategy achieves a speedup between $9\%$ and $64\%$ compared to the cells-batched approach in the discontinuous Galerkin case. With increasing polynomial degrees, the operator becomes more compute intensive, diminishing the advantage.

For cubic shape functions, the difference to hypercube elements with tensor product shape functions is best visible. They profit from the sum-factorization algorithm, decreasing the computational complexity. Furthermore, the increased number of DoFs per cell reduces overheads.

\subsection{Comparison to Matrix-Based Finite Elements and Curvilinear Elements} \label{sec:matrix_based}
With the basic performance properties of matrix-free methods identified, we now compare them against matrix-based approaches. This comparison is crucial for understanding the computational efficiency of each method. By examining their performance characteristics, we can identify scenarios where one method may offer advantages over the other. 

\begin{figure}
    \centering
   \begin{subfigure}{0.5\textwidth}
  \centering
  \includegraphics[width=0.99\linewidth]{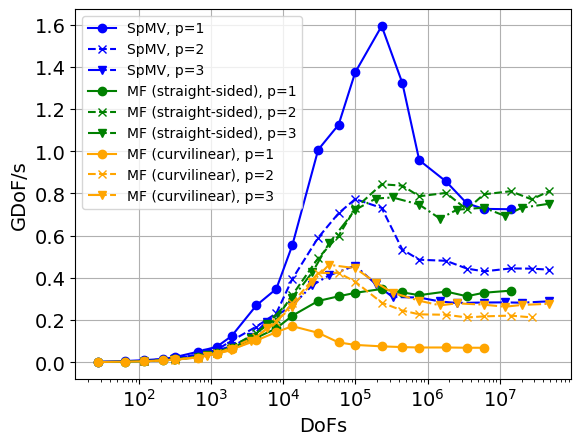}
  \caption{Continuous Galerkin (CG)}
  \label{fig:throughput_matrix_based_CG}
\end{subfigure}%
\begin{subfigure}{.5\textwidth}
  \centering
  \includegraphics[width=0.99\linewidth]{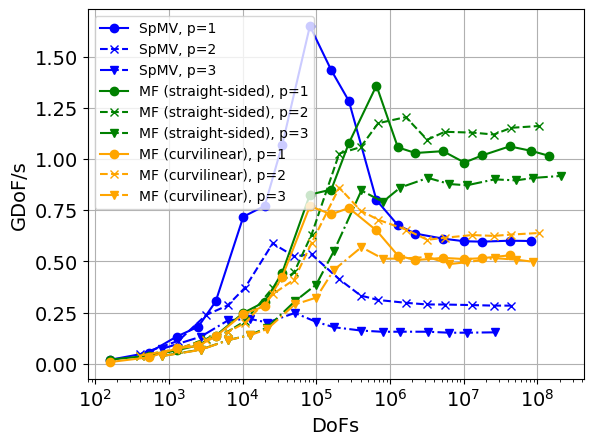}
  \caption{Discontinuous Galerkin (DG)}
  \label{fig:throughput_matrix_based_DG}
\end{subfigure}

\vskip\baselineskip

\begin{subfigure}{0.5\textwidth}
  \centering
  \includegraphics[width=0.99\linewidth]{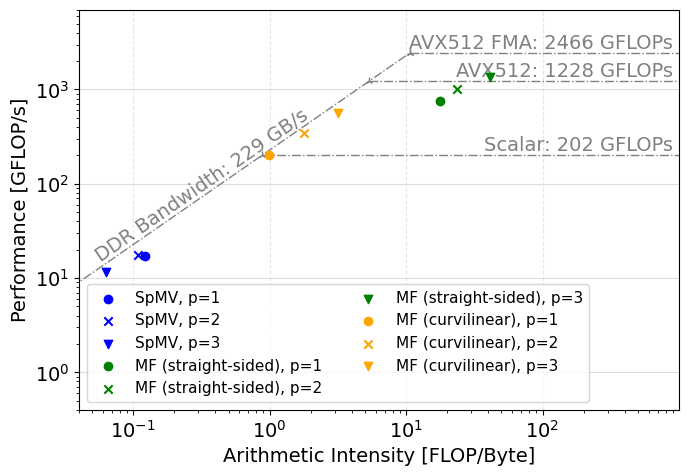}
  \caption{Continuous Galerkin (CG)}
  \label{fig:roofline_matrix_based_curelinear_CG}
\end{subfigure}%
\begin{subfigure}{.5\textwidth}
  \centering
  \includegraphics[width=0.99\linewidth]{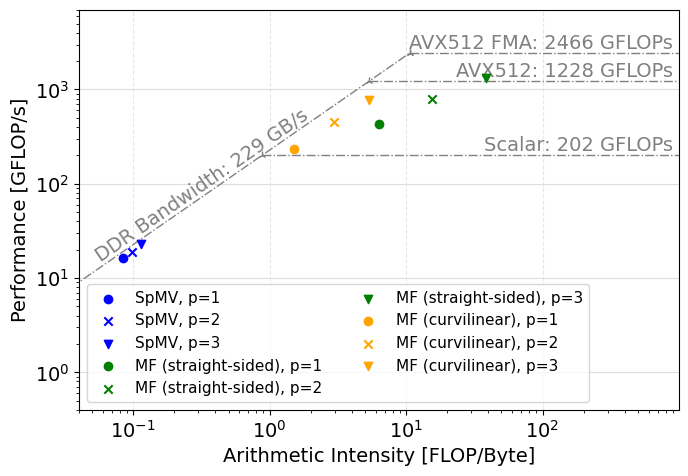}
  \caption{Discontinuous Galerkin (DG)}
  \label{fig:roofline_matrix_based_curelinear_DG}
\end{subfigure}
    \caption{Throughput comparison (top) and roofline model (bottom) of matrix-based implementation and matrix-free implementation on straight-sided and curvilinear elements (scalar Poisson problem). The throughput studies were performed until the system exhausted its available memory, and the roofline plot data were taken from the largest problem size executed prior to memory exhaustion.}
    \label{fig:matrix_based_curelinear}
\end{figure}

The matrix-based code is implemented using \software{Trilinos}~\cite{trilinos}. Figure~\ref{fig:matrix_based_curelinear} compares the throughput of the matrix-based methods with that of the matrix-free methods for a single operator evaluation. For linear continuous elements, the matrix-based method outperforms the matrix-free code, independent of the number of DoFs. At $p=1$, the matrix-free method generates overhead due to repeated work on shared degrees of freedom and has a relatively high memory access per DoF, see Figure~\ref{fig:arith_mem}. At higher polynomial degrees, this overhead is masked by the increased computational intensity in the matrix-free methods, whereas sparse matrices suffer from the denser coupling between degrees of freedom. Pronounced cache effects are observed for the matrix-based implementations with a high throughput at around $10^5$ DoFs, whereas large sizes are strongly limited by main memory bandwidth.

The matrix-free method achieves a higher throughput for the discontinuous elements across all polynomial degrees in the saturated regime. The main reason is that the assembled matrix is less sparse compared to the continuous case. The roofline model in Figure~\ref{fig:roofline_matrix_based_curelinear_DG} illustrates the various performance limits.

For the curvilinear element shapes, a separate Jacobian matrix must be loaded for every quadrature point. The data transfer per cell (Equation~\eqref{eq:mem_trans_cg}) gets extended by the term $9 \cdot n_{q} \cdot 8$, the Jacobian matrix of size $3 \times 3$, with 8 bytes per entry in the case of double precision, for every quadrature point. 
This increases the data access by factors of~32,~17, and~13 for polynomial degrees $p=1,2, \text{ and }3$, respectively. As shown in Figure~\ref{fig:matrix_based_curelinear}, this shifts the arithmetic intensity and leads to a memory-bound scenario, accentuating cache effects. The same applies to the discontinuous case.
In Figure~\ref{fig:throughput_matrix_based_CG} at~$p=2$, the matrix-based code demonstrates higher throughput in the saturated regime compared to the curvilinear matrix-free method. Conversely, for the discontinuous case at polynomial degrees greater than one, the matrix-free code outperforms the matrix-based approach; at~$p=1$ the matrix-free curvilinear elements stay competitive compared to the matrix-based version. 

\begin{figure}
    \centering
   \begin{subfigure}{0.5\textwidth}
  \centering
  \includegraphics[width=0.99\linewidth]{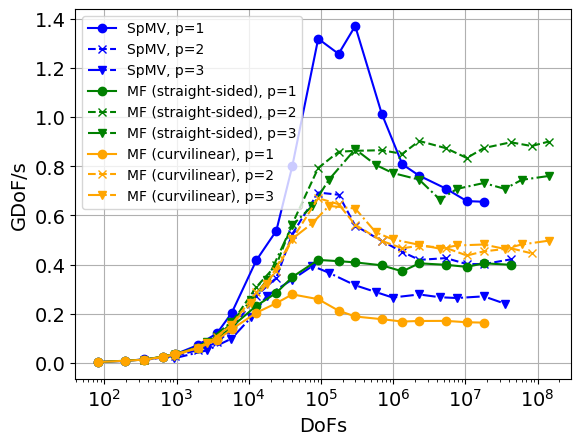}
  \caption{Continuous Galerkin (CG)}
  \label{fig:throughput_components_batched_CG}
\end{subfigure}%
\begin{subfigure}{.5\textwidth}
  \centering
  \includegraphics[width=0.99\linewidth]{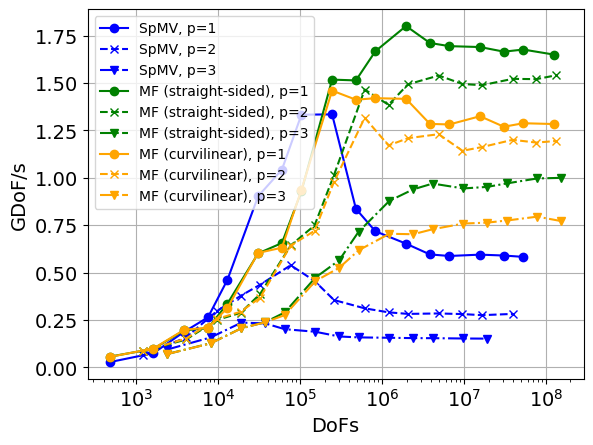}
  \caption{Discontinuous Galerkin (DG)}
  \label{fig:throughput_components_batched_DG}
\end{subfigure}

\vskip\baselineskip
\begin{subfigure}{0.5\textwidth}
  \centering
  \includegraphics[width=0.99\linewidth]{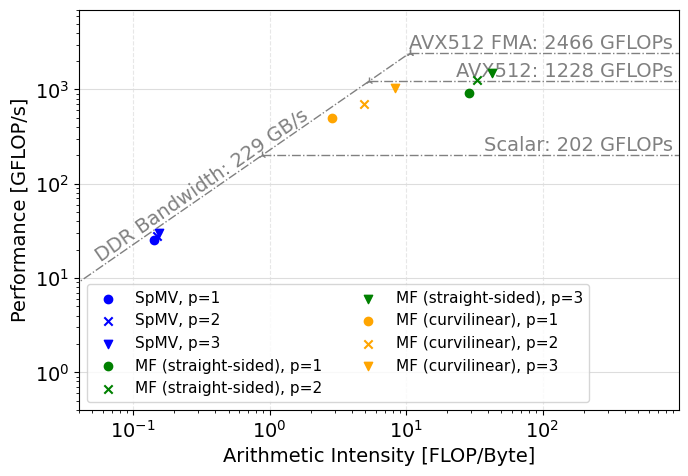}
  \caption{Continuous Galerkin (CG)}
  \label{fig:roofline_components_batched_CG}
\end{subfigure}%
\begin{subfigure}{.5\textwidth}
  \centering
  \includegraphics[width=0.99\linewidth]{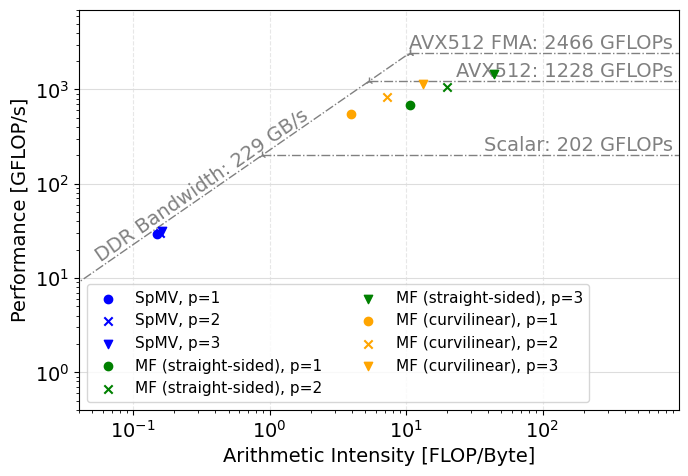}
  \caption{Discontinuous Galerkin (DG)}
  \label{fig:roofline_components_batched_DG}
\end{subfigure}
    \caption{Throughput (top) and roofline analysis (bottom) of matrix-based implementation and matrix-free implementation on straight-sided and curvilinear elements. The problem is a vector-valued Poisson problem with~$d=3$ components.}
    \label{fig:matrix_based_curelinear_components_batched}
\end{figure}

Figure~\ref{fig:matrix_based_curelinear_components_batched} illustrates a similar trend for a vector-valued Poisson equation with three components, optimized using the components-batched strategy. The matrix-based method outperforms the matrix-free code for linear continuous straight-sided and curvilinear elements, see Figure \ref{fig:throughput_components_batched_CG}. However, at higher polynomial degrees, the matrix-based code struggles as more non-zero entries are introduced into the sparse matrix. The same applies to discontinuous elements.
In both cases, a higher throughput of the matrix-free operator is recorded.

In the roofline model (Figures \ref{fig:roofline_components_batched_CG} and \ref{fig:roofline_components_batched_DG}), the arithmetic intensity of matrix-free algorithms increases compared to the scalar Poisson operator, as metric data can be reused for all components.

\begin{remark}
The throughput of the matrix-based implementation depends on the number of non-zero entries in the sparse matrix. We analyzed the behavior of a vector-valued Poisson operator with no coupling between the different components. The behavior changes if coupling between the components is considered. For a problem with coupling between components, e.g.~by a symmetric gradient $\nabla u + (\nabla u)^\text{T}$, the number of non-zero entries per DoF in a sparse matrix increases by a factor of 3 (assuming three components in a 3D setting), which leads to lower throughput compared to the matrix-free operator for every considered case of continuous and discontinuous straight-sided and curvilinear elements,  whereas for the matrix-free operator the coupling merely leads to more operations at the quadrature points (see Equations \eqref{eq:matrix_free_cg} and \eqref{eq:matrix_free_dg}). For further details, see, e.g.,\cite{Munch2024}.
\end{remark}

\subsection{Modified Integration for the Poisson Operator}
The Poisson equation, being a second-order elliptic variational problem, permits the use of an integration rule that is exact for polynomials of degree~$2p-2$ while maintaining its convergence rate~\cite{Ciarlet1976}. This allows for the selection of integration rules with fewer quadrature points in the continuous case.

However, this does not apply to the face integrals for discontinuous Galerkin discretizations. As the jump operators in the penalty term do not involve derivatives, a quadrature rule that is exact up to order $2p$ is needed. Applying a reduced quadrature on the faces decreases the convergence rate of iterative solvers. Hence, in the discontinuous case, only cell integrals are treated with the modified integration while keeping the standard integration for face terms.

\begin{figure}[h!]
    \centering
   \begin{subfigure}{0.5\textwidth}
  \centering
  \includegraphics[width=0.99\linewidth]{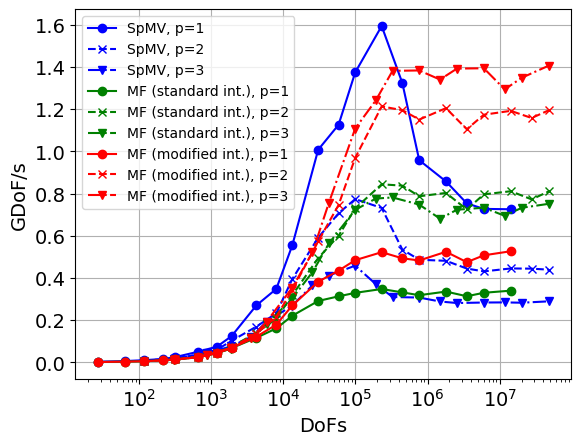}
  \caption{Continuous Galerkin (CG)}
  \label{fig:throughput_underintegration_CG}
\end{subfigure}%
\begin{subfigure}{.5\textwidth}
  \centering
  \includegraphics[width=0.99\linewidth]{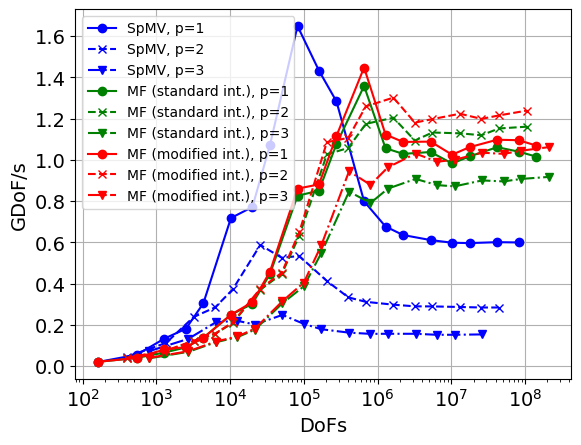}
  \caption{Discontinuous Galerkin (DG)}
  \label{fig:throughput_underintegration_DG}
\end{subfigure}
    \caption{Throughput comparison of matrix-based implementation and matrix-free implementation with standard and modified integration rules (scalar Poisson problem, straight-sided elements).}
    \label{fig:throughput_underintegration}
\end{figure}

Equations~\eqref{eq:arithm_cg} and~\eqref{eq:arithm_dg} indicate that the number of operations depends linearly on the number of integration points. Figure~\ref{fig:throughput_underintegration} illustrates the effect of reducing the number of quadrature points. For the linear continuous operator, which is instruction-bound, reducing the number of quadrature points results in a $64\%$ increase in throughput, yet the matrix-free operator is still outperformed by the matrix-based implementation. In contrast, the operator at $p=3$, which is compute bound, nearly doubles the performance from~$0.75 \text{ GDoF/s}$ to~$1.4\text{ GDoF/s}$. 
For discontinuous elements, the impact of reduced integration points is less pronounced, as the operations on the faces make up around $75\%$ of the arithmetic operations.

\begin{remark}
    The Jacobian matrices of curvilinear elements are high-order polynomials~\cite{Kirby2003, Mengaldo2015}. As no exact integration with the modified integration rule is possible, this case is not further considered in this work.
\end{remark}

\begin{remark}
    In this chapter, only the throughput of the operator action for matrix-based and matrix-free realization was considered. Another factor is the memory consumption of the methods and the time it takes to actually assemble a sparse matrix. For large problems, the matrix-free operator is more efficient in terms of memory use, as the global matrix does not need to be held in memory. The differences can be seen in Figure~\ref{fig:throughput_components_batched_DG}, where up to a factor of $\times 10$ more DoFs can be handled by the same system without running out of main memory.
\end{remark}

\section{Applications}
\label{sec:application}
After optimizing the matrix-vector product, the method is integrated into a multigrid framework to solve the Poisson equation efficiently. In a second step, we demonstrate the characteristics of the operator in the applied problem of airflow in the human lung.

\subsection{Poisson Equation with Hybrid Multigrid Preconditioning}
\label{sec:poisson}
The discontinuous Galerkin method builds a natural ansatz space for the pressure Poisson equation, which is relevant for solving Navier--Stokes problems. Therefore, we consider a scalar Poisson problem on the domain~$\Omega = [-1,1]^3$ with the solution
\begin{equation}
    u(\textbf{x}) = \sin(3 \pi x_1) \sin(3 \pi x_2) \sin(3 \pi x_3) \text{.}
\end{equation}
The right-hand side and the boundary conditions are derived from the solution. To generate the coarse grid, the cube is subdivided into five tetrahedral elements. The implementation is available in the open-source software \software{ExaDG}~\cite{exadg}, which is based on \software{deal.II}. A preconditioned conjugate gradient solver is employed, utilizing a hybrid multigrid strategy as described in~\cite{Fehn2019}.
The hybrid MG preconditioner incorporates polynomial multigrid ($p$-MG), which decreases the polynomial degree by one on each MG-level and geometric multigrid ($h$-MG). The $h$-MG levels are constructed using the grid refinement levels, employing shortest-interior-edge refinement as described in~\cite{Zhang1988}. For discontinuous elements, the problem is additionally transferred to an auxiliary continuous space ($c$-MG)~\cite{Antonietti2017}. The transfer operator between the MG levels is likewise implemented in a matrix-free fashion~\cite{Munch2022}. The coarse grid problem is solved with one V-cycle of the algebraic multigrid (AMG) solver from the~\software{ML}-package~\cite{ml-guide}, with a Chebyshev smoother.

In each preconditioning step, a V-cycle is performed using a Jacobi preconditioned Chebyshev smoother~\cite{Adams2003,Sundar2014}, which demonstrates favorable properties in a matrix-free context~\cite{Munch2022,Kronbichler2018}. On each level five pre--smoothing steps and five post--smoothing steps are conducted. The solver runs in mixed-precision~\cite{Kronbichler2012,Kronbichler2019b}, where the multigrid V-cycle runs in single precision to increase throughput while all other operations in the conjugate gradient solver run in double precision.  

While the solver primarily works with the matrix-free operator, it switches to the matrix-based operator when it offers higher throughput (see Section \ref{sec:matrix_based}). Additionally, the AMG solver requires the assembled system matrix. The matrix size is reduced by the $c$-MG, $p$-MG and $h$-MG transfers relative to the fine-level problem.

As demonstrated in Tables~\ref{tab:n_E_10_DG_MG} and~\ref{tab:n_E_10_CG_MG}, the number of iterations remains constant under mesh refinement. Thereby, the fractional iteration number $n_{10}$ is defined as the number of iterations needed to reduce the residual by ten orders of magnitude, see~\cite{Fehn2019}. Denoting by sequences of $\{c, p, h\}$ the order in which different transfer options are combined from fine to coarser levels, Table~\ref{tab:n_E_10_DG_MG} examines the discontinuous case for the $cph$-MG, $chp$-MG, and $phc$-MG preconditioners. Results for the $hpc$-MG preconditioner are not shown as they are nearly identical to the $phc$-MG scheme. Although all multigrid strategies (i.e., $c$-MG, $h$-MG, $p$-MG, and any combination) exhibit constant iteration counts under mesh refinement, strategies that perform the $c$-transfer at the finest level achieve about one-half to one-third of the iteration numbers compared to strategies beginning with $p$-transfer or $h$-transfer, consistent with results for hypercube elements~\cite{Fehn2019}. The $chp$-MG and $cph$-MG preconditioners reduce the residual by approximately an order of magnitude per cycle. Using a pure $c$-MG strategy provides around~$25\%$ of the throughput compared to $cph$-MG or $chp$-MG, due to the increased number of unknowns in the coarse-grid problem handled by the AMG solver.

\begin{table}
\centering
\caption{Fractional iteration numbers~$n_{10}$ and throughput~$E_{10}$ of the DG MG-solver, needed to reduce the residual by ten orders of magnitude for a scalar Poisson problem on a cube mesh subdivided with tetrahedral elements under different refinement levels.}
\begin{tabular}{|c|r|c|c|c|c|c|c|}
\hline
$p$ & \multicolumn{1}{c|}{$n_{\mathrm{DoF}}$} &  \multicolumn{3}{c|}{$n_{10} \left[ - \right]$} &  \multicolumn{3}{c|}{$E_{10} \left[10^5 \; \mathrm{DoF}/\mathrm{s}/\mathrm{core}\right]$} \\
\cline{3-8}
 & & $cph$-MG &$chp$-MG & $phc$-MG  & $cph$-MG & $chp$-MG & $phc$-MG \\
\hline
1 & 81\,920 & 10.0 & 10.0 & 21.3      & 0.97 & 0.97 & 0.45  \\
1 & 655\,360 & 10.0 & 10.0 & 23.5      & 2.26 & 2.26 & 0.89  \\
1 & 5\,242\,880 & 10.2 & 10.2 & 25.0     & 2.26 & 2.26 & 0.86  \\
1 & 10\,240\,000 & 10.3 & 10.3 & 24.8     & 2.09 & 2.1 & 0.84 \\
\hline
2 & 204\,800 & 9.6 & 9.6 & 22.8   &1.51 & 1.53 & 0.56   \\
2 & 1\,638\,400 & 9.7 & 9.6 & 24.1     &2.70 & 2.86 & 0.84  \\
2 & 13\,107\,200 & 9.8 & 9.5 & 25.3    &2.49 & 2.68 & 0.74  \\
2 & 25\,600\,000 & 9.8 & 9.5 & 25.3  & 2.37 & 2.53 & 0.72  \\
\hline
3 & 409\,600 & 8.9 & 8.9 & 20.9        & 1.57 & 1.57 & 0.52  \\
3 & 3\,276\,800 & 8.8 & 8.7 & 22.7          & 2.27 & 2.40 & 0.75 \\
3 & 26\,214\,400 & 8.7 & 8.6 & 24.2        & 2.26 & 2.43 & 0.63 \\
3 & 51\,200\,000 & 8.7 & 8.7 & 23.9       & 2.20 & 2.35 & 0.63  \\
\hline
\end{tabular}

\label{tab:n_E_10_DG_MG}
\end{table}

\begin{figure}
    \centering

\begin{subfigure}{.49\textwidth}
  \centering
  \includegraphics[width=0.99\linewidth]{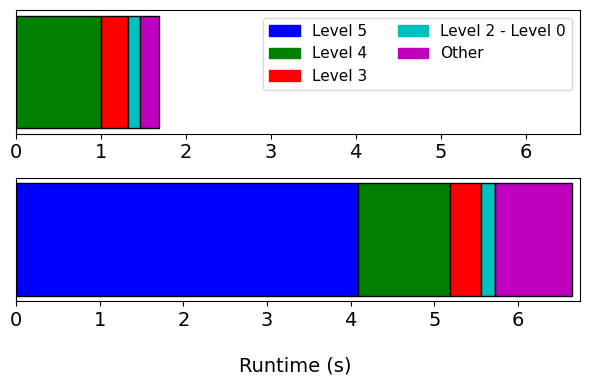}
  \caption{Runtime of the $ph$-MG version (continuous elements, top) and $cph$-MG version (discontinuous elements, bottom) at $p=3$ on a grid with 2.56 million cells.}
  \label{fig:poisson_runtime_DG_CG}
\end{subfigure}
   \begin{subfigure}{0.49\textwidth}
  \centering
  \includegraphics[width=0.99\linewidth]{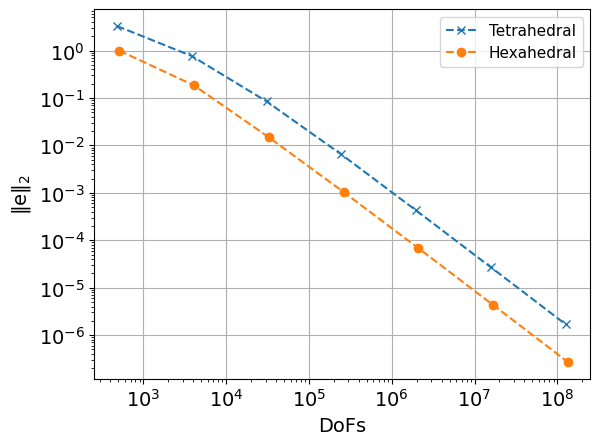}
  \caption{Relative $L^2$ error at $p=3$ under mesh refinement.}
  \label{fig:relative_error_poisson}
\end{subfigure}
    \caption{Left: Runtime of the multigrid solver on the CG and DG problem. Right: Relative $L^2$ discretization error for tetrahedral compared to hexahedral mesh.}
    \label{fig:poisson_MG_results}
\end{figure}

Table ~\ref{tab:n_E_10_CG_MG} presents the same problem setups using continuous elements on identical meshes, yielding similar throughput. The difference in absolute time to solution in Figure~\ref{fig:poisson_runtime_DG_CG} is attributable to the number of DoFs. While at $p$=3 the continuous discretization has around 4.6 unique DoFs per cell, the discontinuous discretization results in 20 DoFs per cell. 

\begin{table}
\centering
\caption{Fractional iteration numbers~$n_{10}$ and throughput~$E_{10}$ of the CG MG-solver, needed to reduce the residual by ten orders of magnitude for a scalar Poisson problem on a cube mesh subdivided with tetrahedral elements under different refinement levels.}
\begin{tabular}{|c|r|c|c|c|c|}
\hline
$p$ & \multicolumn{1}{c|}{$n_{\mathrm{DoF}}$} &  \multicolumn{2}{c|}{$n_{10} \left[ - \right]$} &  \multicolumn{2}{c|}{$E_{10} \left[10^5 \; \mathrm{DoF}/\mathrm{s}/\mathrm{core}\right]$} \\
\cline{3-6}
 & & $ph$-MG & $hp$-MG & $ph$-MG & $hp$-MG \\
\hline
1  & 4\,241   & 5.3 & 5.3 & 0.23 & 0.23 \\
1  & 30\,497  & 5.6 & 5.6 & 0.57 & 0.58 \\
1  & 230\,977 & 6.5 & 6.5 & 0.92 & 0.91 \\
1  & 446\,441 & 6.2 & 6.2 & 0.97 & 0.97 \\
\hline
2  &   30\,497 & 6.3 &  5.6 & 0.83 & 0.89 \\
2  &  230\,977 & 7.1 &  6.0 & 1.41 & 1.85 \\
2  & 1\,797\,249 & 7.3 &  6.1 & 1.67 & 2.57 \\
2 & 3\,491\,281  & 7.3 &  6.1 & 1.85 & 2.67 \\
\hline
3 &         99\,249  &  6.1  & 7.0 & 1.46 & 1.24 \\
3 &        765\,281  &  6.8  & 7.1 & 2.07 & 2.21 \\
3 &       6\,009\,537  &  7.1  & 7.2 & 2.13 & 2.48 \\
3 &      11\,694\,521  &  7.0  & 7.4 &  2.07 & 2.33 \\
\hline
\end{tabular}

\label{tab:n_E_10_CG_MG}
\end{table}

Figure~\ref{fig:poisson_runtime_DG_CG} provides a detailed breakdown of the runtime contributions from various multigrid levels, including the coarse grid solver on level 0 and the time spent in the remaining operations of the conjugate gradient solver (‘Other’). Each multigrid level involves~10 operator evaluations:~4 for pre-smoothing,~5 for post-smoothing with the Chebyshev iteration, and one to compute the residual before restriction. The operator evaluation contributes $79\%$ of the runtime on the different levels, resulting in less time spent on the coarse compared to the fine levels and showing the importance of optimizing the matrix-vector product. The computation of the residual dominates the remaining time spent in the conjugate gradient solver.

As depicted in Figure~\ref{fig:relative_error_poisson}, the convergence rate remains consistent with the problem on a hexahedral grid. The offset between the two graphs is attributed to the worse aspect ratios of the initial tetrahedral grid, which leads to worse error constants. 

\subsection{Incompressible Navier--Stokes Equations in a Human Lung Geometry}
\begin{figure}
    \centering
     \begin{subfigure}{0.33\textwidth}
  \centering
    \includegraphics[width=0.99\textwidth]{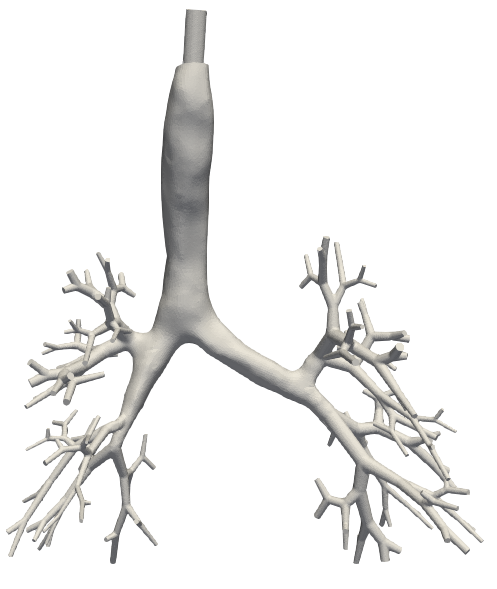}
    \caption{Lung model.}
    \label{fig:lung}
    \end{subfigure}%
 \begin{subfigure}{0.33\textwidth}
  \centering
  \includegraphics[width=0.99\linewidth]{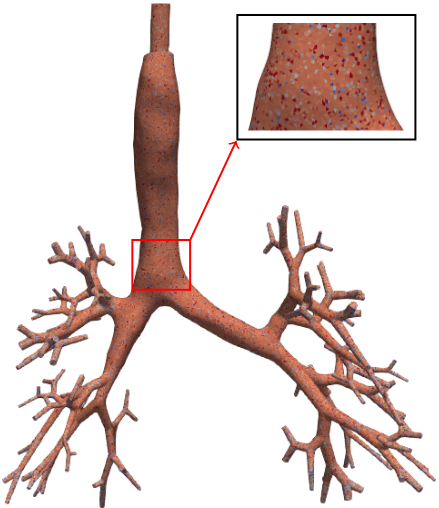}
  \caption{Unordered mesh.}
  \label{fig:lung_unordered_complete}
\end{subfigure}%
\begin{subfigure}{.33\textwidth}
  \centering
  \includegraphics[width=0.99\linewidth]{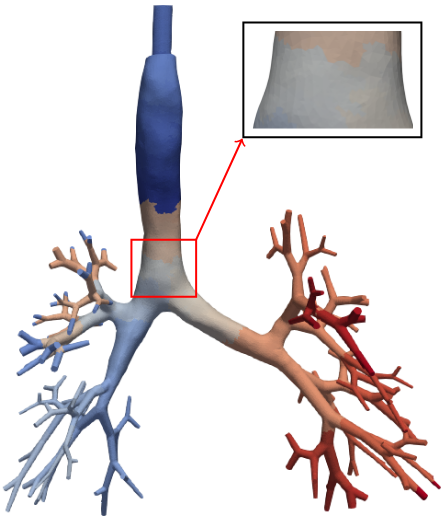}
  \caption{Reordered mesh.}
  \label{fig:lung_reordered_complete}
\end{subfigure}
\caption{Patient specific model of a lung of a preterm infant (left). The mesh as generated by a grid generator (middle) and the hierarchically reordered version (right). The colors in the two figures on the right represent the traversal order of the cells, with blue indicating the first and red the last in the sequence.}
\label{fig:lungs}
\end{figure}

As DG methods are attractive solvers for incompressible turbulent flows~\cite{Chidyagwai2010,Krank2016, Fehn2018b}, as an application example, we look at the flow through a patient-specific model of the human lung of a preterm infant, seen in Figure~\ref{fig:lung}~\cite{Roth2016}. Given the increased geometric complexity, we first examine its impact on the operator.

\subsubsection{Throughput on Unstructured Grids}
To demonstrate the impact of the hierarchical mesh reordering, we compare the throughput across three different mesh configurations: a structured tetrahedral mesh on a cube geometry (as discussed in~\Cref{sec:thoughput}), the unstructured lung mesh as generated by a grid generator (`unordered') and the same lung mesh after hierarchical reordering according to Chapter~\ref{sec:reordering} (`reordered'), see Figures~\ref{fig:lung_unordered_complete} and~\ref{fig:lung_reordered_complete}. The grid has 5.1 million cells and no grid refinement is used. To distribute the grid across MPI processes, \software{METIS} \cite{METIS} is employed.

\begin{figure}
    \centering
    \includegraphics[width=0.7\textwidth]{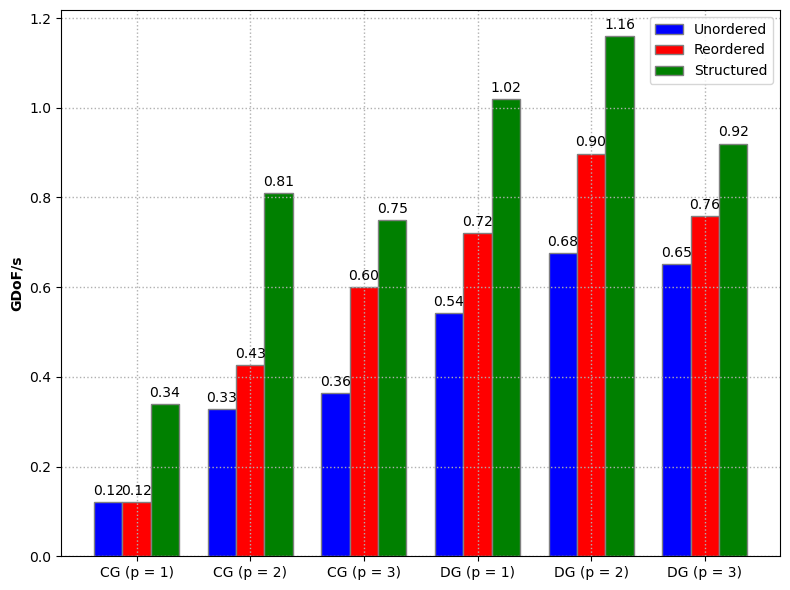}
    \caption{Throughput on structured, unstructured and an unstructured but hierarchical reordered grid of a scalar Poisson operator.}
    \label{fig:throughput_unstructured}
\end{figure}

Figure~\ref{fig:throughput_unstructured} illustrates the throughput of the Poisson operator on these different meshes. The unstructured lung mesh exhibits poor cell ordering, resulting in scattered data access and reduced caching efficiency, seen in Figure~\ref{fig:lung_unordered_complete} as random-like color mixes. The performance difference arises from the higher memory traffic on the unstructured grid, as detailed in Table~\ref{tab:original_reordered_grid}. Hierarchical reordering mitigates this issue by enhancing data locality, increasing cache reuse and reducing bandwidth pressure. The performance gap between the structured and unstructured mesh is expected, as compression of mapping data is more effective on the cube geometry. With increasing polynomial degrees, the difference decreases as the operator tends to be compute-bound. The performance gap between the structured and reordered meshes is smaller for the components-batched strategy compared to the cells-batched strategy, as more data on a cell is reused.

\begin{table} 
\begin{center}
\caption{Measured memory transfer on the unstructured and hierarchically reordered grids in byte/DoF for the cells-batched strategy.}
\begin{tabular}{ |c|c|c|c|c| }
 \hline
   degree & CG unordered & CG reordered &DG unordered & DG reordered \\ 
   \hline
1  &  826  &  846 &  203  &  146 \\
2  &  220  &  162  &  133  &  76 \\
3  & 208 &  77  &  102  &  55 \\
 \hline
\end{tabular} \label{tab:original_reordered_grid}
\end{center}
\end{table}

\subsubsection{Scaling of the Poisson and Helmholtz Operator} \label{sec:scaling}

The Navier--Stokes equations considered here are given by
\begin{equation}
   \frac{\partial \mathbf{u}}{\partial t} + (\mathbf{u} \cdot \nabla) \mathbf{u} = -\nabla p + \nu \nabla^2 \mathbf{u} + \mathbf{f}\text{,}
\end{equation}
with the incompressibility constraint of~$\nabla \cdot \mathbf{u} = 0$. It is solved on the lung geometry with the high-order dual splitting scheme from~\cite{Karniadakis1991}:

\begin{align}
  \frac{\gamma_0 \mathbf{u}^* - \sum_{i = 0}^{J-1} (\alpha_i \mathbf{u}^{n-i})}{\Delta t} &= -\sum_{i = 0}^{J-1} \beta_i (\mathbf{u}^{n-i} \cdot \nabla) \mathbf{u}^{n-i} +  \mathbf{f}^{n+1} \label{eq:dual_splitting_1}\\
   -\nabla^2 p^{n+1} &= -\frac{\gamma_0}{\Delta t} \nabla \cdot \mathbf{u}^* \label{eq:dual_splitting_pressure_poisson}\\ 
  \mathbf{u}^{**} &= \mathbf{u}^* - \frac{\Delta t}{\gamma_0} \nabla p^{n+1}   \label{eq:dual_splitting_2}\\
 \frac{\gamma_0}{\Delta t} \mathbf{u}^{n+1} - \nu \nabla^2 \mathbf{u}^{n+1}  &= \frac{\gamma_0}{\Delta t} \mathbf{u}^{**}  \label{eq:dual_splitting_helmholtz}
\end{align}
where $\textbf{u}^*$ is a extrapolated velocity and $\textbf{u}^{**}$ is a projected velocity.

\begin{figure}
    \centering
   \begin{subfigure}{0.5\textwidth}
  \centering
  \includegraphics[width=0.99\linewidth]{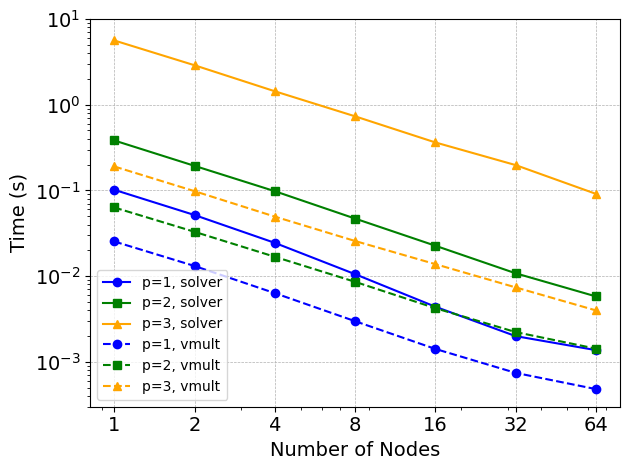}
  \caption{Scaling of the Helmholtz solver.}
  \label{fig:time_vs_n_procs}
\end{subfigure}%
\begin{subfigure}{.5\textwidth}
  \centering
  \includegraphics[width=0.99\linewidth]{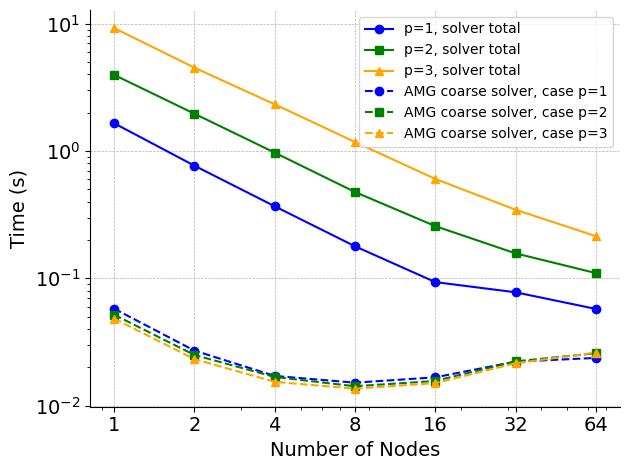}
  \caption{Scaling of the Poisson solver.}
  \label{fig:time_vs_n_procs_poisson}
\end{subfigure}
    \caption[]{Left: Scaling behavior of the vector-valued Helmholtz solver with 2,3 and 18 conjugate gradient iterations at $p$=1, 2 and 3 and a single operator application (vmult). Right: Scaling of the Poisson solver and of the AMG coarse grid solver.}    \label{fig:strong_scaling_helmholz}
\end{figure}

In the preceding sections, we have demonstrated effective strategies for preconditioning and solving the pressure Poisson Equation~\eqref{eq:dual_splitting_pressure_poisson}. Building upon this foundation, we now turn to the computation of a single timestep of the Navier--Stokes operator. This involves additionally computing Equations~\eqref{eq:dual_splitting_1} and~\eqref{eq:dual_splitting_2} and solving the vector-valued Helmholtz like Equation \eqref{eq:dual_splitting_helmholtz}. To this end, we use a discontinuous Galerkin discretization and solve the resulting system using a conjugate gradient solver preconditioned with a point--Jacobi method. To accelerate the solution process, the components-batched strategy is employed in order to solve for the velocity-field. The right side of Equation~\eqref{eq:dual_splitting_helmholtz} is set to unity, with the viscosity of air assumed to be $10^{-5} \frac{\text{m}^2}{\text{s}}$ and a timestep of $10^{-6}\text{s}$ is utilized. We show the scaling of the solver as a limit case for small timesteps. The experiments are done on Intel Xeon Platinum 8360Y CPUs, see Table \ref{tab:xeon_gold_platinum}.
Figure \ref{fig:strong_scaling_helmholz} shows perfect scaling up to 32 nodes at $p=1$ and up to 64 nodes for $p=2$ and $p=3$ for the Helmholtz operator and the conjugate gradient solver. Shown are 2, 3, and 18 iterations of the conjugate gradient solver for $p=$1, 2, and 3, respectively. In the linear case, each process operates on 13300 DoFs (with 64 nodes in total), approaching the latency limit of inter-process communication within the conjugate gradient solver~\cite{Kronbichler2023}. Likewise, the matrix-vector product shows similar behavior, as well as the operations done in Equations~\eqref{eq:dual_splitting_1} and~\eqref{eq:dual_splitting_2}, for which a similar throughput is achieved.

However, the pressure Poisson Equation~\eqref{eq:dual_splitting_pressure_poisson}, solved using the hybrid MG preconditioner, establishes the limit for strong scaling \cite{Krank2016, Fehn2018}. For the $cp$-MG preconditioned pressure Poisson solver, the scaling is good up to 16 and 32 nodes, respectively. Limiting is the AMG coarse grid solver, which takes between $0.03 s$ and $0.155 s$ accumulated over all 14, 13, and 12 solver calls at $p=$1, 2, and 3, respectively, which falls in line with the results of previous studies, see~\cite{Kronbichler2021}. The coarse grid problem consists of around one million DoFs, limiting the effectiveness of the distributed computation. No speedup of the AMG solver can be seen upwards of eight nodes; further optimizations of the AMG solver are beyond the scope of this work.

\section{Conclusions} \label{sec:conclusions}
\begin{sloppypar}
In this study, we have presented a node-level performance optimization approach for continuous and discontinuous Galerkin methods on unstructured tetrahedral grids. Our primary innovation lies in optimizing the operator evaluation, which includes several key techniques such as leveraging explicit data parallelism and Single Instruction, Multiple Data capabilities. Two key strategies were proposed to enhance computational efficiency: batching the degrees of freedom across multiple cells or components into a matrix to enable efficient dense matrix-matrix multiplications with custom compute kernels. Our results show that operators with cubic shape functions are compute-bound, reaching over $60\%$ of the peak performance. The gap to the performance of pure matrix-matrix products is due to separate read/write and compute phases during operator evaluation. Operators with linear shape functions are instruction scheduling bound. Compared to operators with higher-order polynomials, they perform fewer computations due to the lower number of degrees of freedom per element, such the operator cannot compensate as much for the read/write phases. The matrix-free operator achieves a higher throughput for nearly all test cases compared to the matrix-based scheme, except for scalar linear continuous elements and the case of curvilinear elements in specific scenarios. Additionally, hierarchical grid reordering enhances data locality for accessing neighboring cells, thereby improving throughput. Further, we demonstrate that a modified integration rule increases the throughput of the continuous operator, though this benefit did not extend to its discontinuous counterpart.

We have shown the effectiveness of our approach through a series of numerical experiments. The matrix-free approach integrates seamlessly with a geometric multigrid preconditioner, making it a robust tool to solve the (pressure) Poisson equation. The Helmholtz operator demonstrated good strong scalability, providing a robust solution technique for large-scale fluid simulations on unstructured tetrahedral grids. 

Future work will focus on extending these techniques to prismatic and pyramid elements for mixed-mesh support. We anticipate that for higher-order polynomials, sum-factorization techniques on simplex elements are a more powerful method for operator evaluations, falling in line with techniques used on hexahedral grids. Additionally, extending the implementation to GPU architectures could yield further performance improvements, as the parallel architectures excel in the types of operations central to this study. Yet we expect the same constraints identified in Chapter~\ref{sec:scaling} to limit performance in solving large-scale PDE problems.
\end{sloppypar}

\section*{Acknowledgments}\begin{sloppypar}
The research presented in this paper was partly funded by the German Ministry of Research, Technology and Space through project ``PDExa: Optimized software methods for solving partial differential equations on exascale supercomputers'', grant agreement no.~16ME0637K, and the European Union -- NextGenerationEU and by BREATHE, an ERC-2020-ADG Project, Grant Agreement ID 101021526.

The authors gratefully acknowledge contributions by Peter Munch, Richard Schussnig, the \software{deal.II} community and collaborations within the project PDExa. The authors gratefully acknowledge access to HPC resources provided by the Erlangen National High Performance Computing Center (NHR@FAU) of the Friedrich-Alexander-
Universit\"at Erlangen-N\"urnberg (FAU). NHR funding is provided by federal and Bavarian state
authorities. NHR@FAU hardware is partially funded by the German Research Foundation (DFG),
440719683.
\end{sloppypar}

\bibliographystyle{abbrv}
\bibliography{references}
\end{document}